\documentclass[12pt]{article}
\usepackage{amsmath,amsfonts,latexsym,amsthm,amssymb}
\newcommand{\D}{\mathbb D_t^{(\alpha )}}
\newcommand{\Rn}{\mathbb R^n}
\newcommand{\PP}{\mathbf P}
\newcommand{\CM}{\mathbb C^m}
\numberwithin{equation}{section}

\DeclareMathOperator{\const}{const}
\DeclareMathOperator{\R}{Re}
\DeclareMathOperator{\I}{Im}
\begin{document}
\newtheorem{prop}{Proposition}
\newtheorem{lem}{Lemma}
\newtheorem*{teo}{Theorem}
\pagestyle{plain}
\title{Fractional-Hyperbolic Systems}
\author{Anatoly N. Kochubei\footnote{This work was supported in part by Grant No. 01-01-12 of the National Academy of Sciences of Ukraine under the program of joint Ukrainian-Russian projects.}
\\ \footnotesize Institute of Mathematics,\\
\footnotesize National Academy of Sciences of Ukraine,\\
\footnotesize Tereshchenkivska 3, Kiev, 01601 Ukraine\\
\footnotesize E-mail: \ kochubei@i.com.ua}
\date{}
\maketitle

\vspace*{3cm}
\begin{abstract}
We describe a class of evolution systems of linear partial differential equations with the Caputo-Dzhrbashyan fractional derivative of order $\alpha \in (0,1)$ in the time variable $t$ and the first order derivatives in spatial variables $x=(x_1,\ldots ,x_n)$, which can be considered as a fractional analogue of the class of hyperbolic systems. For such systems, we construct a fundamental solution of the Cauchy problem having exponential decay outside the fractional light cone $\{(t,x):\ |t^{-\alpha}x|\le 1\}$.
\end{abstract}
\vspace{2cm}


\newpage
\section{Introduction}

Evolution equations with fractional time derivatives are among central objects of the modern theory of partial differential equations -- due both to various physical applications (dynamical processes in fractal and viscoelastic media \cite{Ma,MK,U}) and to the rich mathematical content of this subject; see, for example, the monographs \cite{EIK,KST} and references therein. The theory has reached sufficient maturity to investigate not only specific classes of equations like fractional diffusion and diffusion-wave equations, but to study general systems of fractional-differential equations. In particular, Heibig \cite{H} obtained an $L^2$-existence theorem for general fractional systems while the author \cite{K12} investigated a class of fractional-parabolic systems.

In this paper we turn to a fractional version of the class of first order hyperbolic systems with constant coefficients (for the classical material see \cite{BS,CH,GS,V}). A prototype is the fractional diffusion-wave equation
\begin{equation}
\left( \mathbb D_t^{(\beta )}u\right) (t,x)=\Delta u(t,x),\quad t\in (0,T],x\in \Rn ,
\end{equation}
where $1<\beta <2$, $\mathbb D_t^{(\beta )}$ is the Caputo-Dzhrbashyan fractional derivative, that is
\begin{equation}
\left( \mathbb D_t^{(\beta )}u\right) (t,x)=\frac{1}{\Gamma (2-\beta )}\frac{\partial^2}{\partial t^2}\int\limits_0^t (t-\tau )^{-\beta +1}u(\tau ,x)\,d\tau -t^{-\beta +1}\frac{u_t(0,x)}{\Gamma (2-\beta )}-t^{-\beta }\frac{u(0,x)}{\Gamma (1-\beta )};
\end{equation}
below we explain how to reduce the equation (1.1) to a system of equations of order $\beta /2\in (0,1)$ in $t$ and the first order in spatial variables.

This equation was studied by many authors (\cite{Fu,Ha,LMP,Ma96,Ps,SW} and others). For the Cauchy problem corresponding to the initial conditions
\begin{equation}
u(0,x)=u_0(x),\quad u_t(0,x)=0,
\end{equation}
a solution is obtained as a convolution with the Green kernel,
\begin{equation}
u(t,x)=\int\limits_{\Rn }G(t,x-\xi )u_0(\xi )\,d\xi ,
\end{equation}
where, as it was proved in \cite{Fu} for $n=1$, and in \cite{Ps} for the general case,
$$
|G(t,x)|\le Ct^{-\frac{\beta n}2}\gamma_n(|x|t^{-\beta/2})E(|x|t^{-\beta/2}),
$$
$$
\gamma_n(z)=\begin{cases}
1, & \text{ if $n=1$};\\
|\log z|, & \text{ if $n=2$};\\
z^{-n+2}, & \text{ if $n\ge 3$},\end{cases}
$$
$$
E(z)=\exp \left( -az^{\frac{2}{2-\beta }}\right),\quad C,a>0
$$
(here and below the letters $C,a$ will denote various positive constants).

For a comparison, note that for the wave equation corresponding formally to $\beta =2$, with the initial conditions (1.3), a counterpart of the kernel $G$ from (1.4) is a distribution supported on the light cone $\{|x|t^{-1}\le 1\}$, if $n$ is even, or its boundary, if $n$ is odd. There is a similar behavior, for example, for some symmetric hyperbolic systems of the first order.

Thus in the fractional case, though the fundamental solution of the Cauchy problem (FSCP) is not concentrated on the set $\left\{ |x|t^{-\frac{\beta}2}\le 1\right\}$, it decays exponentially outside it possessing a kind of weak hyperbolicity property. Of course, the equation (1.1) ``interpolates'' between the heat and wave equations, and possesses some ``parabolic'' properties too, as it is clear from (1.4).

In this paper we describe a class of general time-fractional systems of partial differential equations with constant coefficients, of the form
\begin{equation}
\left( \D u_j\right) (t,x)=\sum\limits_{k=1}^m P_{jk}\left( i\frac{\partial}{\partial x}\right) u_k(t,x),\quad j=1,\ldots ,m,\ 0<t\le T,x\in \Rn,
\end{equation}
where $0<\alpha <1$,
$$
\left( \D \varphi\right) (t)=\frac1{\Gamma (1-\alpha )}\left[ \frac{\partial}{\partial t}\int\limits_0^t (t-\tau )^{-\alpha}\varphi (\tau )\,d\tau -t^{-\alpha}\varphi (0)\right],
$$
$P_{jk}$ are first order polynomials with complex coefficients from the derivatives $i\frac{\partial}{\partial x_1},\ldots ,i\frac{\partial}{\partial x_n}$ ($i=\sqrt{-1}$), for which, as in the above example, there exists a FSCP with an exponential decay away from the set $\left\{ |x|t^{-\alpha}\le 1\right\}$. This property implies the existence of solutions for the initial functions with exponential growth. We call such systems {\it "fractional-hyperbolic"}. We follow the techniques by Gel'fand and Shilov \cite{GS} and Friedman \cite{Fr}, with modifications needed to cover the fractional situation.

\section{Auxiliary Results}

{\bf 2.1.} {\it Derivatives of the Mittag-Leffler functions}. The Mittag-Leffler function
\begin{equation}
E_\alpha (z)=\sum\limits_{l=0}^\infty \frac{z^l}{\Gamma (\alpha l+1)}
\end{equation}
is a fractional calculus counterpart of the exponential function. The function $u(t)=E(\lambda t^\alpha)$ is a solution of the Cauchy problem $(\mathbb D^{(\alpha )}u)(t)=\lambda u(t)$, $u(0)=1$ ($0<\alpha <1$, $\lambda \in \mathbb C$); this representation of the solution remains valid for matrix-valued solutions in the case where $\lambda$ is a matrix.

The generalized Mittag-Leffler functions
$$
E_{\alpha ,\gamma}(z)=\sum\limits_{l=0}^\infty \frac{z^l}{\Gamma (\alpha l+\gamma )}
$$
and
$$
E_{\alpha ,\gamma}^\rho (z)=\sum\limits_{l=0}^\infty \frac{(\rho )_l}{\Gamma (\alpha l+\gamma)}\frac{z^l}{l!},
$$
where
$$
(\rho )_l=\rho (\rho +1)\cdots (\rho +l-1),l\ge 1;\quad (\rho )_0=1,
$$
are also very useful.

In this paper we will need estimates for derivatives of the function (2.1). It is known (see the identities (1.8.23) and (1.9.1) in \cite{KST}) that
\begin{equation}
\frac{d^k}{dz^k}E_\alpha (z)=k!E_{\alpha ,1+\alpha k}^{k+1}(z)
\end{equation}
and
$$
E_{\alpha ,\gamma}^\rho (z)=\frac1{\Gamma (\rho )}{}_1\Psi_1\Biggl[ \begin{matrix}
(\rho ,1) \\ (\gamma,\alpha )\end{matrix}\Biggl| z\Biggr] .
$$
Here ${}_1\Psi_1$ is the Wright function,
$$
{}_1\Psi_1\Biggl[ \begin{matrix}
(\rho ,1) \\ (\gamma,\alpha )\end{matrix}\Biggl| z\Biggr] =\sum\limits_{l=0}^\infty \frac{\Gamma (l+\rho )}{\Gamma (\alpha l+\gamma)}\frac{z^l}{l!}.
$$

The asymptotic behavior of ${}_1\Psi_1$, as $z\to \infty$, was studied in \cite{Wr,Br}. The results for our case are as follows. If $|\arg z|\le \frac{\pi \alpha }2-\varepsilon$, $0<\varepsilon <\frac{\pi \alpha }2$, then
\begin{equation}
E_{\alpha ,\gamma}^\rho (z) \sim \const \cdot \exp (z^{1/\alpha }),\quad z\to \infty .
\end{equation}

If $|\arg (-z)|\le (1-\frac{\alpha}2)\pi -\varepsilon$, $0<\varepsilon <(1-\frac{\alpha}2)\pi$, then
\begin{equation}
E_{\alpha ,\gamma}^\rho (z) \sim Q(-z),\quad z\to \infty ,
\end{equation}
where $Q(z)$ is the sum of residues of the function
$$
s\mapsto z^s\Gamma (-s)\Gamma (s+\rho )/\Gamma (\alpha s+\gamma )
$$
at the points $s=-\rho -\nu$, $\nu =0,1,2,\ldots$. This means that
\begin{equation}
\left| E_{\alpha ,\gamma}^\rho (z)\right| \le C|z|^{-\rho },\quad |z|\ge 1, |\arg (-z)|\le (1-\frac{\alpha}2)\pi -\varepsilon .
\end{equation}

For the remaining region, $\frac{\alpha \pi}2-\varepsilon\le \pm \arg z \le \frac{\alpha \pi}2+\varepsilon$, where $\varepsilon >0$ is small enough,
\begin{equation}
E_{\alpha ,\gamma}^\rho (z) \sim Q(\mp z)+\exp (z^{1/\alpha }).
\end{equation}

It follows from (2.2), (2.4), and (2.5) that for any $\varepsilon >0$,
\begin{equation}
\left|\frac{d^k}{dz^k}E_\alpha (z)\right| \le \begin{cases}
Ce^{\R z^{1/\alpha }}, & \text{ if $|\arg z|\le \frac{\alpha \pi}2$},\\
C_\varepsilon (1+|z|)^{-k-1}, & \text{ if $|\arg z|>\frac{\alpha \pi}2+\varepsilon$}.\end{cases}
\end{equation}
By (2.6),
\begin{equation}
\left|\frac{d^k}{dz^k}E_\alpha (z)\right| \le C,\quad \text{if $|\arg z|\ge \frac{\alpha \pi}2$}.
\end{equation}

\bigskip
{\bf 2.2.} {\it Eskin's extension of the Paley-Wiener-Schwartz theorem}. The Paley-Wiener-Schwartz theorem gives an interpretation of a subclass of the class of entire functions of exponential type as the set of Fourier transforms of distributions with compact supports. Eskin \cite{Es} considered the case of entire functions of order greater than 1.

\medskip
\begin{prop}[Eskin]
Let $Q(s)=Q(s_1,\ldots ,s_n)$, $s_1,\ldots ,s_n\in \mathbb C$, be an entire function of order $p>1$ satisfying the inequality
\begin{equation}
|Q(\sigma +i\tau )|\le C(1+|\sigma |)^le^{a|\tau |^p},\quad \sigma, \tau \in \Rn ,
\end{equation}
where $l$ is a nonnegative integer. Consider $Q(\sigma )$, $\sigma \in \Rn$, as a distribution from $\mathcal S'(\Rn )$. Then its Fourier transform $\widetilde{Q}$ has the form
\begin{equation}
\widetilde{Q}(x)=\sum\limits_{k=1}^N R_k\left( \frac{\partial }{\partial x}\right) F_k(x)
\end{equation}
(in the distribution sense) where the orders of differential operators $R_k$ do not exceed $l+n+1$, if this number is even, or $l+n+2$ otherwise, $F_k$ are continuous functions satisfying, for any $\varepsilon >0$, the inequality
\begin{equation}
|F_k(x)|\le C_\varepsilon e^{-(b-\varepsilon )|x|^{p'}},\quad k=1,\ldots ,N,
\end{equation}
$\frac{1}p+\frac{1}{p'}=1$, $b=\frac{1}{p'}\left( \frac{1}{ap}\right)^{\frac{p'}p}$.
\end{prop}

\medskip
In a further generalization we consider the case where the function $Q$ depends on an additional parameter $t\in [0,T]$, is continuous in $(s,t)$, and both the functions $Q$ and $\D Q$ ($0<\alpha <1$) have the upper bounds (2.9) with all the parameters independent of $t$. Then in (2.1), the operators $R_k$ do not depend on $t$, the functions $F_k$ and $\D F_k$ are continuous in $(x,t)$ and have the upper bounds as in (2.11) with all the parameters independent of $t$. This follows easily from the proof given in \cite{Es} and was noticed (for $\alpha =1$) by Friedman (\cite{Fr}, Theorem $3'$).

\section{Fractional Hyperbolicity}

{\bf 3.1.} {\it Definition}. Let us consider a system (1.5). It is convenient to use its vector form
\begin{equation}
\left( \D U\right) (t,x)=\PP \left(i\frac{\partial }{\partial x}\right)U(t,x),\quad 0<t\le T,x\in \Rn,
\end{equation}
where $U$ is a function with values in $\CM$, $\PP$ is a $m\times m$ matrix whose elements are first order differential operators. In the ``dual'' representation, $\PP (s)$ is a matrix whose elements are polynomials in the variables $s_1,\ldots ,s_n$ of degree 1.

Denote by $\lambda_1(s),\ldots ,\lambda_m (s)$ the roots of the characteristic equation
$$
\det (\PP (s)-\lambda I)=0.
$$
Let
\begin{equation}
\Lambda_\alpha (s)=\max\limits_{|\arg \lambda_k(s)|\le \alpha \pi /2}\R \lambda_k^{1/\alpha }(s).
\end{equation}
If $|\arg \lambda_k(s)|>\alpha \pi /2$ for all $k=1,\ldots ,m$, then we set $\Lambda_\alpha (s)=0$. In any case, $\Lambda_\alpha (s)\ge 0$.

We call the system (3.1) {\it fractional-hyperbolic}, if
\begin{equation}
\Lambda_\alpha (s)\le C\left( |\tau |^{1/\alpha }+\log (|\sigma |+1)\right),\quad s\in \mathbb C^n,
\end{equation}
where $s=\sigma +i\tau$, $\sigma ,\tau\in \Rn$.

\bigskip
{\bf 3.2.} {\it Special classes of systems}. Let us consider systems satisfying the condition
\begin{equation}
\R \lambda_k (s)\le a|\tau |+b \quad (k=1,\ldots ,m)
\end{equation}
where $a,b\ge 0$. Then (3.3) is satisfied. Indeed, if $|\arg \lambda_k(s)|\le \alpha \pi /2$, then
$$
\R \lambda_k^{1/\alpha }(s)=\left| \lambda_k(s)\right|^{1/\alpha }\cos (\frac1\alpha \arg \lambda_k(s)).
$$
An elementary investigation shows that $\cos (\frac1\alpha \varphi)\le \left(\cos \varphi\right)^{1/\alpha}$, if $|\varphi|\le \alpha \pi /2$. Therefore
$$
\R \lambda_k^{1/\alpha }(s)\le \left( \R \lambda_k (s)\right)^{1/\alpha }\le (a|\tau |+b)^{1/\alpha}\le C\left( |\tau |^{1/\alpha}+1\right),
$$
which implies (3.3).

Note that the condition (3.4) is slightly stronger than the Gel'fand-Shilov hyperbolicity condition for differential systems \cite{GS}.

An important class of fractional-hyperbolic systems is that of ``symmetric systems''
\begin{equation}
\left( \D U\right) (t,x)+\sum\limits_{\nu =1}^nA_\nu \frac{\partial U(t,x)}{\partial x_\nu}+BU(t,x)=0
\end{equation}
where $A_\nu$, $\nu =1,\ldots ,n$ are Hermitian matrices, $B$ is an arbitrary matrix. In this case
$$
\PP(s)=i\left( \sum\limits_{\nu =1}^ns_\nu A_\nu +iB\right).
$$

Let $\mu_k(s)$, $k=1,\ldots ,n$, be the eigenvalues of the matrix $L(s)=\sum\limits_{\nu =1}^ns_\nu A_\nu +iB$. Then $\R \lambda_k(s)=\I \mu_k(s)$ for each $k$. By Hirsch's theorem (see Theorem 1.3.1 in Chapter 3 of \cite{MM}), $|\mu_k(s)|\le nM$ where $M$ is the maximum of absolute values of elements of the matrix $\dfrac{1}{2i}(L(s)-L(s)^*)=\sum\limits_{\nu =1}^n\tau_\nu A_\nu +i(B-B^*)$. This implies the inequality (3.4).

\bigskip
{\bf 3.3.} {\it Reduction of the diffusion-wave equation to a system}. Let us consider the diffusion-wave equation (1.1) (where $1<\beta <2$) with the initial conditions
\begin{equation}
u(0,x)=u_0(x),\quad u_t(0,x)=u_1(x),
\end{equation}
where $u_0,u_1$ are continuous functions. In this section we deal with classical solutions, that is we assume that $u\in C^2$ in the spatial variables, $u\in C^1$ jointly in all the variables, there exists, for each $t>0$, the Riemann-Liouville fractional derivative
$$
\left( D^\beta_{0+,t}u\right) (t,x)=\frac{1}{\Gamma (2-\beta )}\frac{\partial^2}{\partial t^2}\int\limits_0^t (t-\tau )^{-\beta +1}u(\tau ,x)\,d\tau ,
$$
and the equalities (1.1) and (3.6) are satisfied pointwise.

Denote $v(t,x)=\left( \mathbb D^{(\beta /2)}u\right) (t,x)$. On a function $\varphi \in C^1$,
\begin{equation}
\left( \mathbb D^{(\beta /2)}\varphi \right) (t)=\frac{1}{\Gamma (1-\frac{\beta}2)}\int\limits_0^t (t-\tau )^{-\beta /2}\varphi'(\tau )\,d\tau ,
\end{equation}
so that
$$
\left|\left( \mathbb D^{(\beta /2)}\varphi \right) (t)\right| \le Ct^{1-\frac{\beta}2}\to 0,\quad \text{as $t\to 0$}.
$$
Therefore $v(0,x)=0$, and $\mathbb D^{(\beta /2)}v$ coincides with the Riemann-Liouville fractional derivative
$$
\left( D^{\beta /2}_{0+,t}v\right) (t,x)=\frac{1}{\Gamma (1-\frac{\beta}2)}\frac{\partial}{\partial t}\int\limits_0^t (t-\tau )^{-\beta/2}v(\tau ,x)\,d\tau ,
$$

In order to calculate $D^{\beta /2}_{0+,t}v$, we write, as in (3.7),
$$
v(t,x)=I_{0+,t}^{1-\frac{\beta}2}u_t(t,x)
$$
where
$$
\left( I_{0+,t}^\gamma \varphi \right)(t)=\frac1{\Gamma (\gamma )}\int\limits_0^t (t-\tau )^{\gamma -1}\varphi (\tau )\,d\tau ,\quad \gamma >0,
$$
is the Riemann-Liouville fractional integral. Using the identity $I_{0+}^{\gamma_1}I_{0+}^{\gamma_2}=I_{0+}^{\gamma_1+\gamma_2}$ (\cite{KST}, (2.1.30)) we find that
$$
D^{\beta/2}_{0+}v=D^{\beta/2}_{0+}I_{0+}^{1-\frac{\beta}2}u_t=\frac{\partial}{\partial t}I_{0+}^{1-\frac{\beta}2}I_{0+}^{1-\frac{\beta}2}u_t=\frac{\partial}{\partial t}I_{0+}^{2-\beta}u_t=
\frac{\partial}{\partial t}\left( \mathbb D^{(\beta -1)}u\right),
$$
so that
$$
\left( D^{\beta/2}_{0+,t}v\right) (t,x)=\frac{\partial}{\partial t}\left[ \left( D^{\beta-1}_{0+,t}u\right) (t,x)-\frac{u_0(x)}{\Gamma (2-\beta )}t^{-\beta +1}\right]
=\frac{\partial}{\partial t}\left( D^{\beta-1}_{0+,t}u\right) (t,x)-\frac{u_0(x)}{\Gamma (1-\beta )}t^{-\beta }.
$$
Since $\dfrac{\partial}{\partial t}D^{\beta-1}_{0+,t}=D_{0+,t}^\beta$ (\cite{KST}, (2.1.34)), we obtain from the definition (1.2) of $\mathbb D^{(\beta )}$ that
$$
\left( \mathbb D_t^{(\beta /2)}v\right) (t,x)=\left( \mathbb D_t^{(\beta )}u\right) (t,x)+\frac{u_t(0,x)}{\Gamma (2-\beta )}t^{-\beta +1}.
$$

We have found that the Cauchy problem (1.1), (3.6) is reduced to the inhomogeneous problem
\begin{align*}
\left( \mathbb D_t^{(\beta /2)}v\right) (t,x)& =\Delta u(t,x)+u_1(x)\frac{t^{-\beta +1}}{\Gamma (2-\beta )},\\
\left( \mathbb D_t^{(\beta /2)}u\right) (t,x)& =v(t,x),
\end{align*}
$$
v(0,x)=0,\quad u(0,x)=u_0(x).
$$

If $u_1(x)\equiv 0$, so that we consider the conditions (1.3), we can reduce the problem further, to a system of the form (3.5). Namely, we set
\begin{align*}
v^0(t,x)& =\left( \mathbb D_t^{(\beta /2)}u\right) (t,x),\\
v^j(t,x)& =\frac{\partial u(t,x)}{\partial x_j},\quad j=1,\ldots ,n.
\end{align*}
Denoting $V(t,x)=(v^0(t,x),v^1(t,x),\ldots ,v^n(t,x))$ we obtain the system
$$
\left( \mathbb D_t^{(\beta /2)}V\right) (t,x)+\sum\limits_{j=1}^n A_j\frac{\partial V(t,x)}{\partial x_j}=0
$$
where
$$
A_j=\begin{pmatrix}
0 & 0 &\ldots & 0 & -1 & 0 &\ldots & 0\\
0 & 0 &\ldots & 0 & 0  & 0 &\ldots & 0\\
\hdotsfor{8}\\
-1 & 0 &\ldots & 0 & 0 &0 &\ldots & 0\\
\hdotsfor{8}\\
0 & 0 &\ldots & 0 & 0  & 0 &\ldots & 0\end{pmatrix}
$$
($-1$ is in the ($j+1$)-th place).

\section{Fundamental Solution of the Cauchy Problem}

{\bf 4.1.} {\it Resolvent matrix-function}. Looking for a FSCP satisfying (3.1) with the initial condition $U(0,x)=\delta (x)$, we apply formally the Fourier transform in $x$ getting the Cauchy problem
$$
\left( \D \widetilde{U}\right) (t,s)=\mathbf P(s)\widetilde{U}(t,s),\quad \widetilde{U}(0,s)=I,
$$
where $I$ is the unit matrix. Then
\begin{equation}
\widetilde{U}(t,s)=E_\alpha (t^\alpha \mathbf P(s)).
\end{equation}
The function (4.1) is called {\it the resolvent matrix-function} corresponding to the system (3.1).

Below we will denote the inverse Fourier transform (in the sense specified later) of the function (4.1) by $G(t,x)$ leaving the notation $U(t,x)$ for a vector-valued solution with a general initial vector-function.

\medskip
\begin{prop}
For all $t\in (0,T)$, $x\in \Rn$,
\begin{equation}
\left\| E_\alpha (t^\alpha \mathbf P(s))\right\| \le C\left( 1+t^\alpha |s|\right)^{m-1}e^{t\Lambda_\alpha (s)}.
\end{equation}
\end{prop}

\medskip
{\it Proof}. Let us obtain an estimate for $\left\| E_\alpha (t^\alpha \mathcal P)\right\|$ where $\mathcal P$ is an arbitrary matrix with eigenvalues $\mu_1,\ldots ,\mu_m$. Suppose first that all the eigenvalues are different. Then $E_\alpha (t^\alpha \mathcal P)=f(\mathcal P)$ where $f$ is the Newton interpolation polynomial
$$
f(\mu )=b_1+b_2(\mu -\mu_1)+b_3(\mu -\mu_1)(\mu -\mu_2)+\cdots +b_m(\mu -\mu_1)\cdots (\mu -\mu_{m-1}),
$$
and the coefficients are chosen in such a way that $f(\mu_j)=E_\alpha (t^\alpha \mu_j)$, $j=1,\ldots ,m$; see, for example, \cite{Hi} regarding functions of matrices.

It follows from general estimates given in \cite{GS} (Chapter II, Section 6.1) that
$$
|b_k|\le \max\limits_{\mu \in B_k}\left| \frac{d^{k-1}}{d\mu^{k-1}}E_\alpha (t^\alpha \mu )\right|,\quad k=1,\ldots ,m,
$$
where $B_k$ is the smallest convex polygon containing the points $\mu_1,\ldots ,\mu_k$. Using (2.7) and (2.8) we come to the inequality
$$
|b_k|\le Ct^{\alpha (k-1)}e^{tM},\quad M=\max\limits_{|\arg \mu_j|\le \frac{\alpha \pi}2}\R \mu_j^{1/\alpha}.
$$
Therefore
\begin{multline*}
\left\| E_\alpha (t^\alpha \mathcal P)\right\|=\| f(\mathcal P)\|\le Ce^{tM}\left[ 1+t^\alpha (\|\mathcal P\| +|\mu_1|)+t^{2\alpha }(\|\mathcal P\| +|\mu_1|)(\|\mathcal P\| +|\mu_2|)+\cdots \right] \\
\le Ce^{tM}\left( 1+2t^\alpha \|\mathcal P\|+\cdots +(2t^\alpha )^{(m-1)}\|\mathcal P\|^{m-1}\right),
\end{multline*}
since $|\mu_j|\le \|\mathcal P\|$ for all $j$. Here the constant $C$ does not depend on $\mathcal P$ and its eigenvalues. For continuity reasons, the above inequality holds for any matrix $\mathcal P$, not necessarily with distinct eigenvalues.

Now we set $\mathcal P=\PP (s)$. Since $\|\mathcal P\|^2$ does not exceed the sum of squared absolute values of all the elements, and since polynomials appearing as elements of $\PP (s)$ have degrees $\le 1$, we find that $\|\PP (s)\|\le C(1+|s|)$. Substituting and recalling that in this case $M=\Lambda_\alpha (s)$, we obtain the required inequality (4.2).$\qquad \blacksquare$

\medskip
If our system (3.1) is fractional-hyperbolic, that is it satisfies (3.3), then it follows from (4.2) that
\begin{equation}
\left\| E_\alpha (t^\alpha \mathbf P(s))\right\| \le C(1+|\sigma |)^qe^{at|\tau |^{1/\alpha}}
\end{equation}
with some $q>0$. We call the minimal possible nonnegative integer $q$, for which (4.3) holds, {\it the exponent} of the system (3.1). In particular, if our system satisfies (3.4), then $q\le m-1$.

Note that
$$
\D E_\alpha (t^\alpha \mathbf P(s))=\PP (s)E_\alpha (t^\alpha \mathbf P(s)),
$$
so that
$$
\left\| \D E_\alpha (t^\alpha \mathbf P(s))\right\| \le C(1+|\sigma |)^{q+1}e^{at|\tau |^{1/\alpha}},
$$
since the matrix $\PP \left( i\frac{\partial}{\partial x}\right)$ contains only first order differential operators.

\bigskip
{\bf 4.2.} {\it The Cauchy problem}. Let us consider the system (3.1) satisfying (3.3), with the initial condition $U(0,x)=U_0(x)$, $x\in \Rn$. A vector-function $U(t,x)$, $0\le t\le T$, $x\in \Rn$, is called {\it a classical solution} of the Cauchy problem, if it is continuous in $(t,x)$, as well as its first derivatives in $x$, the fractional integral
$$
\left( I_{0+}^{1-\alpha}U\right) (t,x)=\frac1{\Gamma (1-\alpha)}\int\limits_0^t(t-\tau )^{-\alpha}U(\tau ,x)\,d\tau
$$
has the first derivative in $t$, continuous in $(t,x)$, and the equation (3.1) and the initial condition are satisfied pointwise.

Though our aim is to construct a classical solution of the Cauchy problem, at the first stage it will be defined as a generalized solution. Given the estimate (4.3), the further reasoning is very similar to that of \cite{Fr} (in particular, to the case of correctly posed systems with positive genus $\mu$); only the exponent $\dfrac{p_0}{p_0-\mu}$, where $p_0$ is the reduced order, is replaced with $\dfrac1{1-\alpha}$. Therefore
here and in the proof of the theorem below, we do not repeat the lengthy calculations from \cite{Fr} referring instead to specific sections of Friedman's work.

Denote by $W_{\frac1{1-\alpha},c}$ ($c>0$) the Frechet space of such functions $\varphi \in C^\infty (\mathbb R)$ that for any $c'<c$,
$$
\left| \varphi^{(j)}(z)\right| \le C_{j,c'}e^{-(1-\alpha)|c'z|^{\frac1{1-\alpha}}},\quad j=0,1,2,\ldots ,
$$
for all $z\in \mathbb R$. Our basic space $\Phi$ of test functions is the direct product of $n$ copies of $W_{\frac1{1-\alpha},c}$. We will need also a larger space $\Phi^0$, the direct product of $n$ copies of $W_{\frac1{1-\alpha},c-\varepsilon}$, $0<\varepsilon <c$.

The Fourier transform $\mathcal F$ maps $W_{\frac1{1-\alpha},c}$ onto the space $W^{\frac1\alpha ,\frac1c}$ of such entire functions $\psi (s)$, $s=\sigma +i\tau$, $\sigma ,\tau\in \mathbb R$, that for any $d>\frac1c$, $s\in \mathbb C$,
$$
\left| s^k\psi (s)\right| \le C_{k,d}\exp \left( \alpha |d\tau |^{\frac1\alpha}\right),\quad k=0,1,2,\ldots .
$$
Denote $\Psi =\mathcal F\Phi$, $\Psi^0 =\mathcal F\Phi^0$. These spaces have an obvious direct product structure.

It follows from the estimate (4.3) and the description of multipliers in $W^{\frac1\alpha ,\frac1c}$ (\cite{GS}, Section I.2.4) that multiplication by $E_\alpha (t^\alpha \mathbf P(s))$ ($0\le t\le T$) sends $\Psi$ into $\Psi^0$, if $T$ is small enough. Therefore the above multiplication defines an action between the conjugate spaces $\left( \Psi^0\right)'\to \Psi'$. Thus, if $U_0\in \left( \Phi^0\right)'$, then $\widetilde{U_0}=\mathcal FU_0\in \left( \Psi^0\right)'\subset \Psi'$, so that $\widetilde{U}(t,s)=E_\alpha (t^\alpha \mathbf P(s))\widetilde{U_0}(s)$ is a generalized solution over $\Psi$. Therefore
\begin{equation}
U(t,x)=\left( G(t,\cdot )*U_0(\cdot )\right) (x)
\end{equation}
is a generalized solution in the sense of the distribution space $\Phi'$, defined for $t\le T$, if $T$ is small enough.

\medskip
\begin{teo}
There exists such a positive number $\gamma$ that for any initial function $U_0$ possessing continuous derivatives $D^\nu U_0$, $|\nu|\le q+n+3$ ($q$ is the exponent of the system), such that
\begin{equation}
\left| D^\nu U_0(x)\right| \le Ce^{\gamma |x|^{\frac1{1-\alpha}}},\quad x\in \Rn,
\end{equation}
the generalized solution (4.4) is a classical solution of the Cauchy problem, with the estimate
\begin{equation}
|U(t,x)| \le Ce^{\gamma' |x|^{\frac1{1-\alpha}}},\quad 0\le t\le T,x\in \Rn,
\end{equation}
where $C,\gamma'>0$ do not depend on $t,x$. The fundamental solution $G(t,x)$ has the form
\begin{equation}
G(t,x)=\sum_{k=1}^M R_k \left( \frac{\partial}{\partial x}\right) f_k(t,x)
\end{equation}
(the differential operators $R_k$ of orders $\le q+n+3$ are understood in the sense of $\mathcal S'(\Rn )$), $f_k$ are continuous functions satisfying the estimates
\begin{equation}
|f_k(t,x)|\le Ce^{-\gamma_1|t^{-\alpha} x|^{\frac1{1-\alpha}}},\quad \gamma_1>0,
\end{equation}
with the constants independent of $t,x$.
\end{teo}

\medskip
The properties (4.7)-(4.8) express the ``fractional-hyperbolic'' behavior of the class of systems considered in this paper.

\medskip
{\it Scheme of proof}. The representation (4.7) with the estimates (4.8) follow from (4.3) and Proposition 1 (including its generalization described in Section 2.2). Note that the inequality (2.11) implies the estimate
$$
|f_k(t,x)|\le C_\varepsilon e^{-(\gamma_2t^{-\frac{\alpha}{1-\alpha}}-\varepsilon ) |x|^{\frac1{1-\alpha}}},\quad \gamma_2>0,
$$
Taking $\varepsilon <\frac{\gamma_2}2 T^{-\frac{\alpha}{1-\alpha}}$ we come to (4.8) with $\gamma_1=\gamma_2/2$.

It is shown in \cite{Fr}, pages 355-356, that if the convolution $U(t,x)$ defined as a distribution from $\Phi'$ is in fact a continuous function in $(t,x)$ and has, with all the derivatives appearing in (3.1), the upper bounds like that in (4.6) with coefficients $\gamma'<(1-\alpha )c^{\frac1{1-\alpha}}$ where $c$ is the constant involved in the definition of the space $\Phi$, then $U$ is the classical solution.

Next (\cite{Fr}, page 356), in the situation considered here, the convolutions in the distribution sense coincide with the classical ones, so that the estimates of convolutions for functions of exponential growth and exponential decay (\cite{Fr}, page 362) are applicable, and differential operators in expressions like (2.10) can be switched upon smooth convolutors (\cite{Fr}, page 357). Taking into account the assumption (4.5) and using the above convolution estimates we obtain (4.6) for any given $T$, if $\gamma$ is small enough. $\qquad \blacksquare$

\end{document}